\documentclass[12pt]{amsart}
\usepackage{amsfonts}

\theoremstyle{definition}

\theoremstyle{remark}

\newcommand{\const}{\mathop{\rm const}\limits}

\begin{document}

\begin{center}

{\bf WEIGHT HARDY-LITTLEWOOD INEQUALITIES FOR DIFFERENT POWERS.} \par

\vspace{3mm}
{\bf E. Ostrovsky}\\

e - mail: galo@list.ru \\

\vspace{3mm}

{\bf L. Sirota}\\

e - mail: sirota@zahav.net.il \\

\vspace{4mm}

 Abstract. \\
{\it In this short article we obtain the non-asymptotic upper and low estimations for linear and bilinear weight Riesz's functional through the Lebesgue spaces.} \\

\end{center}

\vspace{3mm}


\vspace{3mm}


\vspace{3mm}

\section{Introduction. Statement of problem.}

\vspace{3mm}

 The linear integral operator $ I_{\alpha,\beta,\lambda} f(x), $ or, more precisely, the {\it family of operators} of a view

 $$
u(x) = I_{\alpha,\beta,\lambda}f(x) = I[f](x) = |x|^{-\beta}
\int_{R^d} \frac{f(y) \ |y|^{-\alpha} \ dy}{|x - y|^{\lambda}}
 $$
is called {\it Weight Riesz's integral operator,} or simply {\it Weight Riesz's} potential, or {\it weight fractional integral.} \par
Here $ |x|, \ x \in R^d $ denotes usually Euclidean norm of the vector $ x, \
d = 1,2,3,\ldots; \ \alpha,\beta,\lambda = \const, \alpha,\beta \ge 0, \
\lambda > 0; \ \alpha + \beta + \lambda < d. $ \par

 The bilinear {\it Weight Riesz's functional} $ B_{\alpha,\beta,\lambda}(f,g) = B(f,g) $ may be defined as follows:

 $$
 B(f,g) = \int_{R^d} \int_{R^d} \frac{f(y) \ |x|^{-\beta} \ g(x) \ dx \ dy}
{|y|^{\alpha} \ |x - y|^{\lambda}}.
 $$
 It is evident that

 $$
 B_{\alpha,\beta,\lambda}(f,g) = (I_{\alpha,\beta,\lambda}f, g),
 $$
where $ (f,g)$ denotes ordinary the inner product of a (measurable) functions $ f $ and $ g: $

$$
(f,g) = \int_{R^d} f(x) \ g(x) \ dx,
$$
which is defined, e.g. when $ f \in L_m, \ g \in L_l, \ m,l > 1, 1/m + 1/l = 1,
L_m = L_m(R^d) $ is the classical Lebesgue space of all the measurable functions
$ f: R^d \to R $ with finite norm
 $$
 |f|_m \stackrel{def}{=} \left( \int_{R^d} |f(x)|^m \ dx \right)^{1/m}; \ f \in L_m \ \Leftrightarrow |f|_m < \infty.
 $$
  Obviously,

  $$
  \sup_{|g|_l =1} |B_{\alpha,\beta,\lambda}(f,g)| = |I_{\alpha,\beta,\lambda}f|_m.
  $$

 The operators $ I_{\alpha,\beta,\lambda} $ and correspondingly the functionals
$ B_{\alpha,\beta,\lambda} $
 are used in the theory of Fourier transform, theory of Partial Differential
 Equations, probability theory (study of potential functions for Markovian processes and spectral densities for stationary random fields), in the
 functional analysis, in particular, in the theory of interpolation of operators etc., see for instance \cite{Adams1}, \cite{Bennet1}, \cite{Stein0}, \cite{Stein1}, \cite{Lieb0},\cite{Lieb2}, \cite{Rubin1}, \cite{Opic1}, \cite{Ostrovsky2}, \cite{Ostrovsky3}, \cite{Ostrovsky4}, \cite{Hunt1}, \cite{Krein1} etc. \par

 We denote also $ L(a,b) = \cap_{p \in (a,b)} L_p. $ \par

{\bf We will investigate the estimations of a view:}

$$
|I_{\alpha,\beta,\lambda}f|_q \le K_{\alpha,\beta, \lambda}(p) \ |f|_p,
\eqno(1.a)
$$

$$
|B_{\alpha,\beta,\lambda}(f,g)| \le K_{\alpha,\beta,\lambda}(p) \ |f|_p \ |g|_{q/(q-1)}, \eqno(1.b)
$$
{\bf with asymptotically exact values of coefficient $ K. $ } \par
 Note that the case $ \alpha=\beta=0, $ i.e. the case of classical Riesz potential, is considered in many publications \cite{Adams1}, \cite{Hardy1}, \cite{Lieb0}, \cite{Lieb1}, \cite{Lieb2}, \cite{Lieb3}, \cite{Leoni1}, \cite{Stein1}, \cite{Stein2}, \cite{Talenti1} etc.
 General case $ \alpha^2 + \beta^2 > 0 $ is partially considered, e.g. in \cite{Criado1},\cite{Hardy1}, \cite{Stein0}, \cite{Lieb0}, \cite{Okikiolu1}, \cite{Lacey1}, \cite{Procdoliti1}, \cite{Pick1}. \par
 In the recent publication \cite{Napoli1} is considered the weight Riesz potential for the so-called radial function, i.e. for the functions
$ f(x) $ which dependent only on the Euclidean norm of a vector $ x: f(x) =
H(|x|). $ It is obtained in \cite{Napoli1} the upper bound for the weight
Riesz potential without the constants $ K_{\alpha,beta,\lambda}(p) $ estimations. \par
 We intend to improve the results of these works and to obtain the {\it low} bounds for the functionals $ |B_{\alpha,\beta,\lambda}(f,g)| $ and
 $ |I_{\alpha,\beta,\lambda}(f)|.$ \par

 In order to formulate the main result, we must introduce some notations. We define first of all the following function $ q=q(p) $ as follows:

$$
1+\frac{1}{q} = \frac{1}{p} + \frac{\alpha+\beta+\gamma}{d}. \eqno(2)
$$

We will denote the set of all such a values $ (p,q) $ as $G(\alpha,\beta,\lambda) $ or
for simplicity $ G = G(\alpha,\beta,\lambda). $ \par

{\sc Further we will suppose that }
$ (p,q) \in G(\alpha,\beta,\lambda) = G. $ \par
 It is known \cite{Lieb0}, \cite{Stein0}, \cite{Okikiolu1}
 that the inequalities (1a) and (1b) are possible only in the case when
 $ (p,q) \in G(\alpha,\beta,\lambda). $ \par

 We denote also
$$
p_-:=\frac{d}{d-\alpha}, \ p_+:= \frac{d}{d-\alpha-\lambda};
$$
and correspondingly
$$
q_-:= \frac{d}{\beta+\lambda}, \ q_+ := \frac{d}{\beta},
$$
where in the case $ \beta=0 \ \Rightarrow q_+:= + \infty; $
$$
\kappa = \kappa(\alpha,\beta,\lambda) := (\alpha + \beta + \lambda)/d.
$$
 It is known \cite{Stein0}, \cite{Okikiolu1} that if $ p \in [1,p_-] \cup [p_+,\infty), $ then $ K_{\alpha,\beta, \lambda}(p) = +\infty. $ {\sc Thus, we confine the values $ p $ inside the open interval } $ p \in (p_-,p_+). $ \par
Define the {\it exact value} of the constant $ K_{\alpha,\beta,\lambda}(p), $ i.e. the value

 $$
V(p)= V_{\alpha,\beta\lambda}(p)= \sup_{f \in L(p_-,p_=), f \ne 0}
 \frac{|I_{\alpha,\beta,\lambda}(f)|_q}{|f|_p}. \eqno(3)
 $$

\vspace{3mm}

 We use symbols $C(X,Y),$ $C(p,q;\psi),$ etc., to denote positive
constants along with parameters they depend on, or at least
dependence on which is essential in our study. To distinguish
between two different constants depending on the same parameters
we will additionally enumerate them, like $C_1(X,Y)$ and
$C_2(X,Y).$ The relation $ g(\cdot) \asymp h(\cdot), \ p \in (A,B), $
where $ g = g(p), \ h = h(p), \ g,h: (A,B) \to R_+, $
denotes as usually

$$
0< \inf_{p\in (A,B)} h(p)/g(p) \le \sup_{p \in(A,B)}h(p)/g(p)<\infty.
$$
The symbol $ \sim $ will denote usual equivalence in the limit
sense.\par
We will denote as ordinary the indicator function
$$
I(x \in A) = 1, x \in A, \ I(x \in A) = 0, x \notin A;
$$
here $ A $ is a measurable set.\par
\hfill $\Box$

\bigskip

\section{ Main Result: upper and low bounds for weight Riesz's potential.}

\vspace{3mm}

{\bf Theorem 1.} For the values $ (p,q) \in G(\alpha,\beta,\lambda) $ and
$ p \in (p_-,p_+) $ there holds:
$$
\frac{C_1(\alpha,\beta,\lambda)}{\left[(p-p_-) \ (p_+ - p)\right]^{\kappa} } \le
V_{\alpha,\beta,\lambda}(p) \le \frac{C_2(\alpha,\beta,\lambda)}
{\left[(p-p_-) \ (p_+ - p)\right]^{\kappa} }. \eqno(4)
$$
\vspace{3mm}
Recall that for the values $ (p,q) \notin G(\alpha,\beta,\lambda) $ or for the values
$ p \notin (p_-,p_+) \ V_{\alpha,\beta,\lambda}(p) = \infty. $ \par
{\bf Proof of the upper bound} it follows immediately from \cite{Okikiolu1}, pp. 215-219 after simple calculations.\par
 Another method used in \cite{Ostrovsky3}, \cite{Ostrovsky4} based on the theory of maximal operators, see in \cite{Adams1}, \cite{Lacey1}, \cite{Criado1}, \cite{Pick1}, \cite{Procdoliti1}. \par

{\bf Proof of the low bound.} We will consider two examples of a functions from the set $ L(p_-,p_+).$ \par
\vspace{3mm}
{\bf First example.}
$$
f_0(x) = |x|^{-(d-\alpha)} \ I(|x| > 1).
$$
We find by direct calculations using the multidimensional polar coordinates:

$$
|f_0|_p \asymp \ (p-p_-)^{-1/p} \asymp (p-p_-)^{-(1-\alpha/d)},
p \in (p_-, p_+);
$$
$ |x| \ge 1 \ \Rightarrow $
$$
u_0(x):= I_{\alpha,\beta,\lambda}f_0(x)\ge C^{(1)}_{\alpha,\beta,\lambda} \
|x|^{1-\beta - \lambda} \ |\log |x|| \ I(|x| > 1);
$$

$$
|u_0|_q \ge C^{(2)}_{\alpha,\beta,\lambda} \cdot \times
\left[q - q_0 \right]^{-1-1/q},
$$
recall that $ q = q(p). $ \par
\vspace{3mm}
{\bf Second example.} We put:

$$
g_0(x) = |x|^{-(d-\alpha-\lambda)} \ I(|x| < 1),
$$
and find:

$$
|g_0|_p \asymp (p_+ -p)^{-1/p} \asymp (p_+ -p)^{-(d-\alpha-\lambda)/d};
$$

$$
v_0(x) := I_{\alpha,\beta,\lambda}g_0(x) \ge C^{(3)}(\alpha,\beta,\lambda) \cdot
|x|^{-\beta} \cdot |\log |x|| \cdot I(|x| < 1);
$$

$$
|v_0|_q \ge C^{(4)}(\alpha,\beta,\lambda) \cdot
 (q_+ - q)^{-1-1/q}, \ p \in (p_-, p_+).
$$
\vspace{3mm}
{\bf Third example. Summing.} \par
We define:

$$
h(x):= f_0(x) + g_0(x); \ w(x):= I_{\alpha,\beta,\lambda}[h](x).
$$

 It follows for the function $ h = h(x), $ which belongs to the space $ L(p_-, p_+), $ that for the values $ p $ from the considered interval $ p \in (p_-, p_+) $

$$
\inf_{ p \in (p_-, p_+)} \frac{|w|_q \ [(p - p_-)(p_+ - p)]^{\kappa}}{|h|_p} =:
C^{(6)}(\alpha,\beta,\gamma;d) > 0.
$$

  This completes the proof of theorem 1.\par

\bigskip

\hfill $\Box$


\vspace{5mm}

\begin{thebibliography}{99}

\vspace{4mm}

\bibitem{Adams1}
D.R.Adams, L.I. Hedberg. {\it Function Spaces and Potential Theory.}
Springer Verlag, Berlin, Heidelberg, New York, 1996.
\bibitem{Bennet1}
C. Bennet and R. Sharpley, {\it Interpolation of operators.}
Orlando, Academic Press Inc., 1988.
\bibitem{Carro1}
M. Carro and J. Martin, {\it Extrapolation theory for the real
interpolation method.} Collect. Math. {\bf 33}(2002), 163--186.
\bibitem{Criado1}
A.Criado. {\it On the Lack of Dimension free Estimates in $L^p$ for maximal Functions associated to radial Measures.}
arXiv:0907.4326v1 [math.CA] 24 Jul 2009.
\bibitem{Hardy1}
G.H.Hardy, J.E. Litlewood and G.P\'olya. {\it Inequalities.}
Cambridge, (1952).
\bibitem{Hunt1}
R.A.Hunt. {\it Developments Related to the A.E.Convergence of Fourier Series. }
MAA Studies in Harmonic Analysis. V.13, J. Math. Ass., (1998), p. 20-37.
\bibitem{Krein1}
S.G. Krein, Yu. V. Petunin and E.M. Semenov, {\it Interpolation of
Linear operators.} New York, AMS, 1982.
\bibitem{Lacey1}
M.T.Lacey, K.Moen, C.P\'erez and R.H.Torres. {\it Sharp weighted Bounds for fractional Integral Operators.}
arXiv:0905.3839v1 [math.CA] 23 May 2009.
\bibitem{Ledoux1}
M. Ledoux and M. Talagrand. {\it Probability in
 Banach Spaces.} Springer, Berlin, 1991.
\bibitem{Leoni1}
G.Leoni. {\it A First Course in Sobolev Spaces.} Graduate Studies in Mathematics, Am. Math. Soc., (2009), V. 105, New York, Princeton.
\bibitem{Lieb0}
E.H.Lieb, M.Loss. {\it Analysis. } Graduate Studies in Mathematics, V.14,
Second Edition, Am. Math. Soc., (2003), USA, New York, Princeton.
\bibitem{Lieb1}
R.L.Frank and E.H.Lieb. {\it Inversion Positivity and the sharp Hardy-Littlewood-Sobolev Inequality.}
Electronic Publications, arXiv:0904.4275v1 [math.FA] 27 Apr 2009.
\bibitem{Lieb2}
J.Fr\"olich, R.Israel, E.H.Lieb, B.Simon. {\it Phase transitions and reflections positivity. I. General Theory and long range lattice models.} Commun. Math. Phys.,
62 (1978), no 1, 1-34.
\bibitem{Lieb3}
E.H.Lieb. {\it Sharp constants in the Hardy-Littlewood-Sobolev and related inequalities.}
Ann. of Math., (2), 118 (1983), no 2, 349-374.
\bibitem{Mitrinovich1}
D.S.Mitrinovich, J.E. Pecaric and A.M.Fink. {\it Inequalities involving
Functions and their Integrals and Derivatives.} Kluvner Academic Verlag,
(1996), Dorderecht, Boston, London.
\bibitem{Napoli1}
P.L. De Napoli, I.Drelichman and R.G.Duran. {\it Om weighted Inequalities for
fractional Integrals of radial Functions.} arXiv: 0910.5508v1 [math. CA]
 28 Oct 2009.
\bibitem{Okikiolu1}
G.O.Okikiolu. {\it Aspects of the Theory of Bounded Integral Operators in
$ L^p $ Spaces.} Academic Press, (1971), London - New York.
\bibitem{Opic1}
B.Opic and A.Kufher. {\it Hardy-type inequalities.} Pitman Resaerch Notes in Mathematics
Series, 219 (1990), Longman Scientific and Technical, Harlow, Essex, UK, John Willey and Sons Inc,m New York.
\bibitem{Ostrovsky2}
E. Ostrovsky and L.Sirota. {\it Moment Banach spaces: theory and applications.}
HAIT Journal of Science and Engeneering, {\bf C}, Volume 4, Issues 1 - 2,
pp. 233 - 262, (2007).
\bibitem{Ostrovsky3}
E. Ostrovsky, E.Rogovee and L.Sirota. {\it Riesz's and Bessel's Operators in bilateral grand Lebesgue Spaces. }
 arXiv:0907.3321v1 [math.FA] 19 Jul 2009.
\bibitem{Ostrovsky4}
E. Ostrovsky and L.Sirota. {\it Hardy-Littlewood Inequalities for Riesz Potential:
low Bounds Estimations for different Powers.}
arXiv:0909.5663v1 [math.FA] 30 Sep 2009.
\bibitem{Pick1}
L.Pick. {\it Two-weight weak type maximal inequalities in Orlicz classes.} Studia Matematica, 100 {\bf 3}, (1991), P. 207-218.
\bibitem{Procdoliti1}
G.Procdoliti. {\it Weighted Inequalities and pointwise Estimations for the multilinear
 fractional Integral and maximal Operators.}
arXiv:0907.5210v1 [math.AP] 29 Jul 2009.
\bibitem{Rubin1}
B.Rubin . {\it Fractional Integrals and Potentials.} Pitman Monographs and Surveys in
Pure and Applied Math., (82), Longman, (1988), British Library Publisher House, London.
\bibitem{Stein0}
E.M.Stein and G.Weiss. {\it Fractional Integrals on n-dimensional Euclidean Space.}
J. Math. Mech., {\bf 7}, p. 503-514.
\bibitem{Stein1}
E.M.Stein. {\it Singular Integrals and Differentiability Properties of Functions. }
Princeton University Press, Princeton, (1992).
\bibitem{Stein2}
E.M.Stein and J.O. Str\"omberg. {\it Behavior of maximal functions in $ R^n $ for
large } $ n.$ Arkiv f\"or matematik. Volume 21, no 1, May (1983); Published by
Institut Mittag-Lefler; Djursholm, Sweden, p. 259-269.
\bibitem{Stein3}
E.M.Stein. {\it The development of square finctions in the works of A.Zygmund.}
Bull. Amer. Math. Soc., 7, (1982), p. 359-376.
\bibitem{Stein4}
E.M.Stein. {\it Maximal functions: Spherical means.}Proc. Nat. Acad, Sci. USA,
73(1976), p, 2174-2196.
\bibitem{Steim5}.
E.M.Stein. {\it Harmonic Analysis. Real variable Methods, Orthogonality, and
Oscillatory Integrals.} Princeton, New Jersey, Univ. Press, (1993), p. 12 - 33.
\bibitem{Talenti1}
G.Talenti. {\it Best constants in Sobolev inequality.} Ann. di Matem. Pura Ed Appl.,
110 (1976), 353-372.

\end{thebibliography}
\end{document}